**Andrei Rodin** (Ecole Normale Superieure)

# On Categorical Theory-Building: Beyond the Formal

## 1. Introduction: Languages, Foundations and Reification of Concepts

The term "language" is colloquially used in mathematics to refer to a theory, which grasps common features of a large range (or even all) of other mathematical theories and so can serve as a unifying conceptual framework for these theories. "Set-theoretic language" is a case in point. The systematic work of translation of the whole of mathematics into the set-theoretic language has been endeavoured in 20-th century by a group of mathematicians under the collective name of *Nicolas Bourbaki.* Bourbakist mathematics proved successful both in research and higher education (albeit not in the school math education) and is practiced until today at mathematical departments worldwide.

The project of Bourbaki as well as other attempts to do mathematics "set-theoretically" should be definitely distinguished from the project of reduction of mathematics to set theory defended by Quine and some other philosophers. This latter project is based on the claim that all true mathematical propositions are deducible from axioms of Zermelo-Frenkael set theory with Choice (ZFC) or another appropriate system of axioms for sets, so basically all the mathematics *is* set theory. A working mathematician usually sees this claim as an example of philosophical absurdity on a par with Zeno's claim that there is no motion, and Bourbaki never tried to put it forward. Actually Bourbaki used set theory (more precisely their own version of axiomatic theory of sets) for doing mathematics in very much the same way in which classical first-order logic is used for doing axiomatic theories like ZFC. In other words Bourbaki made set theory a "part of their logic", and then developed specific mathematical theories stipulating new axioms expressed in this extended language. According to the philosophical view just mentioned all Bourbaki's proofs are nevertheless translatable into deductions from axioms of set theory. I shall not go for pros and contras about this controversial claim here but put instead these two questions: (1)Which features of sets make set theory a reasonable candidate for foundations of mathematics? and (2)Which features of sets allow set theory to be an effective mathematical language? Although the two questions are mutually related they are not the same and require different answers.

In order to answer the first question remind Hilbert's *Grundlagen der Geometrie* of 1899 which provides the notion of foundation relevant to the question. Hilbert suggests to think of geometrical points and straight lines as of abstract "things" (of two different types) holding certain relations with required formal properties; one is left free to imagine then these things in any way one likes or not imagine them in any particular way at all. However abstract and unspecific might be the notion of thing involved here one cannot avoid making certain assumptions about it. In order to clarify these assumptions one needs an appropriate "theory of things". Set theory proves appropriate for this purpose: sets provide the standard (Tarski) semantic for classical first-order logic and for theories axiomatised with this logic. Hence the idea to use set theory on a par with formal axiomatic method and the claim that mathematics is ultimately "about" sets. Building of axiomatic set theories becomes then a rather tricky business since any such theory involves an infinite regress: in order to build an axiomatic theory of sets one needs to assume some (usually different) notion of set in advance for semantic purposes. This and other relevant problems about set theory and logic have been scrutinised by mathematicians and philosophers throughout 20th century. I shall not explore this vast issue here but I want to stress the intimate link between sets and formal axiomatic method just explained. True, this method allows for building not only theories of sets but also theories of lines and points (plane geometries), of parts and wholes (mereologies) and of whatnot. However some notion of set (or *class*) is anyway required by all such theories. Stronger technical notions of set corresponding to ZFC and other axiomatic set theories come about when this general requirement is further strengthened for specific mathematical needs. So when the notion of foundations is understood in the sense of Hilbert's *Grundlagen* or similarly the choice of set theory (rather than mereology or anything else) as foundations of mathematics is natural.

Second question requires, in my view, a very different answer. There is more than one reason why set theory is helpful for mathematics but perhaps the following one is the most important. Mathematics doesn't deal with "pure" concepts - whatever this might mean - but deals with concepts embodied into mathematical objects. There are different ways of such embodiment or reification of concepts. The most traditional one is *exemplification*. When Euclid wants to prove that all triangles have certain property *P* he always proceeds in the following way. He takes "for example" just one triangle *ABC*, proves that *ABC* has property *P* and then concludes that *all* triangles have this property. The conclusion is justified when the proof doesn't rely on any specific property *S* of *ABC*, which some other triangle *A'B'C'* might not have. This guarantees that the same proof can be applied to *A'B'C'* or any other triangle. Such "reification through instantiation" is also widely used in today's mathematics as well. It seems to be fairly fundamental for mathematics as we know it through its history. However set theory provides another way of reifying concept. In the pre-set-theoretic mathematics people had among available objects particular triangles like *ABC* and particular (natural) numbers like *5* but they had no special mathematical object corresponding to the general concept of triangle or to the general concept of number (over and above the aforementioned particulars). Entities of this latter kind could be believed to exist somewhere on Platonic heavens but certainly not among mathematical objects. However set theory allowed for making such things up and treating them mathematically. Consider set **T** of all triangles (on a given plane) and set **N** of natural numbers. These sets are *extensions* of their corresponding concepts. They are genuine mathematical objects having certain properties (e.g. cardinality) and allowing for certain operations with them. In particular, **N** can be squared (the square of **N** is the set of all ordered number pairs) and factorised by some equivalence relation. Importantly **T** is not a triangle and **N** is not a natural number: the extension of a given concept is *not* an instance of this concept. But like an instance
the extension is an embodiment of a given concept: it is an object one can play with, i.e. make further constructions. Obviously this second way of reification of concepts was not available for mathematicians until the notion of infinite set was approved by G. Cantor and his followers.

One may remind Occam's Razor and wonder why having more objects is an advantage. The answer is this: it is an advantage as far as it leads to new non-trivial mathematics. This is indeed the case as far as set-theoretic mathematics is concerned. Mathematicians can be interested in a conceptual parsimony but hardly in an ontological parsimony. I cannot see any profit mathematics might get by preventing certain concepts from reification. In particular, since mathematical reasoning involves the notion of infinity anyway it is quite appropriate for this science to develop a calculus of infinities rather than keep the notion of infinity somewhere at the limit of the scope of the discipline without a properly mathematical treatment. The history of mathematics of last two centuries provides numerous examples of successful reification of "ideal elements" of different kinds. Think of ideal points in projective geometry, for example.

## 2. Language of Categories and Categorical Foundations

Unlike set theory category theory has been designed as a language to begin with and only later has been proposed as a foundations by Lawvere and his followers. The notion of category has been first explicitly introduced by Eilenberg and MacLane in their (1945) paper as a purely auxiliary device, and until works of Grothendieck and his school in late fifties, which made an essential use of category theory in algebraic geometry, nobody would consider this theory as anything more than a convenient system of notation. In his classical (1971) MacLane writes:

Category theory starts with the observation that many properties of mathematical systems can be unified and simplified by a presentation with diagrams of arrows. (MacLane 1971, p.1)

Such presentation is often possible because most of mathematical concepts come with a corresponding notion of map (otherwise called transformation or morphism) between tokens falling under a given concept. For example, maps between sets are functions, maps between topological spaces are continuous transformations, maps between groups are group homomorphisms. Such maps are composable in the usual way corresponding to the common intuition behind the notion of transformation. The mathematical notion of category makes this common intuition explicit assuming associativity of composition of maps and existence of identity map for each object. Using the "language of arrows" we may think of, say, natural numbers, not just as of "bare set" **N** = {0, 1, 2, 3, ...} but as a category where numbers are provided with succession maps:

0 --> 1 --> 2 --> 3 --> ...

Identifying number *n* with isomorphism classes of sets having exactly *n* elements and considering classes of maps between these sets as morphisms between numbers one gets a richer category comprising $m\exp n$ different morphisms from any number $n$ to any number $m$, in particular $n\exp n$ morphisms to each number $n$ to itself $n!$ of which are isomorphisms. All these notions can be rather easily cooked set-theoretically; in more involved constructions such a set-theoretic regression is also often (but not always) possible but the advantage of using the language of arrows can be much more important. I shall not talk about specific mathematical matters here but remark that category theory like set theory allowed for reification of certain concepts which earlier could not be reified. Consider the concept of set for example. The extension **S** of this concept is the set of all sets. That **S** turns to be a contradictory notion is a part of the problem but not yet the whole problem (notice that the notion of infinite set before Cantor was believed to be contradicory too). Another part of the problem is that **S** doesn't have interesting properties to be studied and apparently doesn't allow for further non-trivial constructions. However as far as all sets are taken together with *all maps* between them the situation changes. The category of sets **Set** does have interesting specific properties distinguishing it from other categories and also allows for non-trivial constructions (like that of topos). This clearly shows that **Set** is a better embodiment of the general concept of set than **S**. The situation is similar in the case of the concept of group, topological space and many others (see paragraph 6 below). So category theory allows for reification of concepts when set theory doesn't.

However useful category theory might be what has been said so far has no bearing on the problem of foundations. One may assume standard set-theoretic foundations and then construe the language of categories upon it. But why not to think (in particular to think about sets) category-theoretically to start with? The first systematic attempts of this kind has been made by Lawvere in his thesis of 1963 and papers of 1964 and 1966 based on this thesis. In these works Lawvere introduced categories using formal axiomatic method just like Zermelo and Fraenkel did this with sets. This amounts to the following: objects and morphisms are taken as primitives objects holding three primitive relations with intended meaning "domain of", "codomain of", and "composition of" plus the identity relation. Categorical objects and categorical morphisms are treated as belonging to the same type since every categorical object is formally identified with its identity morphism. Lawvere himself avoids speaking about objects and relations in this context taking first purely syntactical viewpoint and after listing the appropriate axioms saying :

By a category we of course understand (intuitively) any structure which is an interpretation of the elementary theory of abstract categories ... (1966, p.4)

In his (1964) Lawvere adjoins to his axiomatic category theory a number of additional axioms making an abstract category "into" the category of sets so that

There is essentially only one category which satisfies these ... axioms ... , namely the category S of sets and mappings. (1964, p. 1506)

Remind that in order to use set theory as foundations one needs "abstract" sets rather than "concrete" sets like sets of points, numbers, etc. What these abstract sets are sets of? Cantor's answer is the following: abstract sets are sets of "pure units" ("lauter Einsen"). Another answer has been later given by Zermelo: abstract sets are sets of sets. This latter answer is obviously more economical conceptually. For a similar reason Lawvere put forward a theory of categories of categories (but not just a general theory of "concrete" categories like categories of sets, groups, etc...) and suggested it as foundations of mathematics (Lawvere 1966).
While in these early papers Lawvere sticks to formal axiomatic method and the corresponding notion of foundations in his more recent paper of 2003 this author takes a different approach and opts for a different notion of foundations understood

... in a common-sense way rather than in the speculative way of the Bolzano-Frege-Peano-Russell tradition. (2003 p. ??)

This change of Lawvere's view seems me remarkable. The intimate link between sets and formal axiomatic method stressed in the previous paragraph suggests that sets cannot be replaced in their foundational role by categories or anything else unless one continue to use this method and applies the corresponding notion of foundation. This, in my view, explains why the idea of making categories into foundations finally led Lawvere to the refusal from the formal method. But Lawvere's reference to common sense hardly solve the problem either. To get rid of "speculative foundations" one needs a new method of theory-building. In his (2003) Lawvere doesn't aim at general solution of this problem but gives a concrete example of how categorical foundations may look like. The principle aim of my paper is to describe a general method of theory-building suggested by category theory. I shall call this new method *categorical* and distinguish it from formal method. Formal views on mathematics and science are usually opposed to more traditional views according to which mathematics and sciences always assume certain "substances" like "number" or "magnitude" as their subject-matters. In today's philosophy of mathematics formal view is associated with mathematical structuralism. I say that categorical method goes beyond formal method (and beyond structuralism) in order to stress that my proposal has nothing to do with the traditional substantialism. As we shall see the mathematical notion of category suggests something genuinely new with respect to the traditional concepts of form and substance.

The rest of the paper is organised as follows. In the next section I stress a distinctive feature of formal axiomatic method (seen against traditional axiomatic method), which concerns the notion of interpretation relevant to mathematics. Then I argue that formal axiomatic method doesn't provide an adequate treatment of mathematical interpretation and introduce the notion of categorical method which does this. Then I analyse some logical aspects of categorical method and conclude with general epistemological arguments in its favour.

## 3. Formal Axiomatic Method in The Nutshell

A today's mathematical student can read in various textbooks that formal axiomatic method invented by Hilbert is nothing but a perfectioned version of the traditional axiomatic method known since Euclid. True, Hilbert certainly had Euclid's *Elements* in mind writing his *Grundlagen*, so his method can be rightly seen as a modification of Euclid's. However I don't think that the description of this modification as perfectioning sheds a lot of light on it. To see clearly what is specific for formal axiomatic method as distinguished from more traditional versions of axiomatic method a historical regression seems me helpful. Soon after the publication of Hilbert's *Grundlagen* in 1899 Frege sent Hilbert a letter (the precise date is missing) containing a severe criticism of Hilbert's approach. Frege's had the following traditional understanding of axiomatic method in mind. A given theory starts with axioms, which are truths taken for granted. These non-demonstrable truths are truths about certain objects. The theory proceeds with inferences from the axioms made according rules of inference, which must be also assumed. As a matter of course for any given theory meanings of all terms used in its axioms and further inferences must be unequivocally fixed once and for all. This general epistemological view dating back to Aristotle has been recently called *classical model of science* (Jong&Betti, forthcoming).

Frege pointed to Hilbert that his *Grundlagen* falls short of meeting the requirements just mentioned and in particular the unequivocality requirement. Here is a quotation from Hilbert's reply to Frege where Hilbert explains his new method in the nutshell:

"You say that my concepts, e.g. "point", "between", are not unequivocally fixed <...>. But surely it is self-evident that every theory is merely a framework or schema of concepts together with their necessary relations to one another, and that basic elements can be construed as one pleases. If I think of my points as some system or other of things, e.g. the system of love, of law, or of chimney sweeps <...> and then conceive of all my axioms as relations between these things, then my theorems, e.g. the Pythagorean one, will hold of these things as well. In other words, each and every theory can always be applied to infinitely many systems of basic elements. For one merely has to apply a univocal and reversible one-to-one transformation and stipulate that the axioms for the transformed things be correspondingly similar. Indeed this is frequently applied, for example in the principle of duality, etc." (cit. by Frege 1971, p.13, italic mine).

Since a point is allowed to "be" (or to "be thought of as") a "system of love and chimney sweeps" (or a beer mug according to another popular Hilbert's saying) - and all this within one and the same theory - Frege's notion of axiomatic method is certainly no longer relevant. But let's look for a serious mathematical reason behind Hilbert's colourful rhetoric. In the end of the quotation Hilbert refers to the duality principle in projective geometry. Given a true proposition of this theory, which involves straight lines and points, one may formally exchange terms "line" and "point" and get another true proposition. (This doesn't reduce to the trivial remark that one may call lines "points" and call points "lines" without changing the given theoretical structure since the original proposition remains true as it stands.) This suggests the following idea: a given mathematical object can be occasionally "thought of" or "interpreted" as another mathematical object. In particular in projective geometry one may "think of lines as points and think of points as lines". Such a liberal treatment of mathematical objects is common in today's mathematics. In the end of 19th century it was not yet common but a number of important examples were already around (I elaborate on one such example in the next section). Hilbert's *Grundlagen* provides a justification for this apparently careless conceptual game. The problem Hilbert addressed can be formulated as follows: How to construe a mathematical concept, which can be occasionally "interpreted" as another concept?; How to formulate a theory in which basic concepts are defined only "up to interpretation"? Hilbert's answer is roughly this. One should first conceive of mathematical objects as bare "things" (possibly of different types) standing in certain relations to each other, and then describe these relations stipulating their formal properties as axioms. Any "system of things", which hold relations satisfying the axioms would be a model of the theory.

As Hilbert makes it explicit in the quoted passage he thinks about re-interpretations of theories as "reversible one-to-one transformations", i.e. as isomorphisms. Hilbert's clearness often absent in later expositions of formal method allows one to see its limits: mutual interpretations of mathematical theories are, generally, not reversible. To build a theory "up to isomorphism" is not the same thing as to build a theory "up to interpretation". For interpretations are, generally, not isomorphisms. Let me now demonstrate this fact using a historical example, which was already available to Hilbert.

4. Irreversible Interpretations

Non-Euclidean geometries emerged in mathematics of 19th century as a result of at least two different developments. The first started in Ancient times and culminated with (Bolyai 1832) and (Lobachevsky 1837). These mathematicians like their predecessors tried to prove the 5-th Postulate of Euclid's Elements (the "Axiom of Parallels") by getting a contradiction from its negations but at certain point changed their attitude and came to a conviction that they were exploring a new vast territory rather than approaching the desired dead end. The second line of development, which I associate with the names of Gauss and Riemann, was relatively recent. Gauss had a genuinely new insight on the old "problem of parallels" guessing a link with geometry of curve surfaces. This allowed his pupil Riemann to build a new generalised concept of geometrical space (Riemanean manifold), which still serves us as the best mathematical description of the physical space-time (in General Relativity).

The two lines of development were brought together by Beltrami in his two prominent papers (1868) and (1868-69). Anachronistically speaking, in his *Saggio* of 1868 Beltrami gave a 3D Euclidean model of plane Lobachevskian geometry. More precisely it was only a partical model where finite segments of geodesics of a surface named by Beltrami pseudo-sphere represented straight line segments of Lobachevskian plane. But Beltrami didn't have the notions of formal theory and model in mind. He first thought he discovered what the Lobachevskian plane was *indeed*: he believed it was the pseudo-sphere. However this conclusion was not quite satisfactory even in Beltrami's own eyes. He didn't notice that his model for plane Lobachevskian gometry was only partial (this was first noticed by Helmholz in 1870, see Kline 1972) but he saw that Lobachevskian 3D geometry couldn't be treated in the same way. So he looked for a better solution. He found it after reading Riemann's (1854) and presented in his *Teoria* of 1869: Lobachevskian space is a Riemanean manifold of constant negative curvature. This holds for spaces of any number of dimensions.

The latter solution apparently makes the talk of interpretation no longer necessary. Let's however see how the result of *Saggio* looks from the point of view of *Teoria*: a 2-manifold of constant negative curvature is partially embeddable into 3D Euclidean space (which is another Riemanean manifold). So we have here two manifolds and a map, which can still be tought of as interpretation as suggested by *Saggio*. The point I want to stress is that this map is not an isomorphism, it is not reversible. It restricts to an isomorphism (a part of Lobachevskian plane is isomorphic to a part of Euclidean space) but the whole construction cannot be conceived on this restricted basis alone: the map in question is a map between two spaces (manifolds) but not between their "parts". As a surface in Euclidean space the pseudo-sphere cannot be "carved out" of this space. One may remark that we are talking about a map between spaces (manifolds) but not about an interpretation between theories, and so this example is not quite relevant to the issue discussed in the previous section. But it is obvious that however the notion of theory is construed in this case the situation remains asymmetric: while Lobachevskian plane geometry can (modulo needed reservations) be explained in or "translated into" terms of Euclidean 3D geometry the converse is not the case. Observe that the mere existence of interpretation *f* :**A**-->**B** of theory **A** in terms of another theory **B** and a backward interpretation *g*:**B**-->**A** is not sufficient for considering *f* as reversible: *f* and *g* should "cancel" each other for it. To give a precise definition one needs to stipulate appropriate "identical interpretations" id$_A$, id$_B$ (which leave **A**,**B** correspondingly "as they are") and then require *fg*=id$_A$ and *gf*=id$_B$, where *fg* and *gf* denote composition interpretations (written here in the direct geometrical order). This standard categorical definition of isomorphism often allows for but doesn't require thinking of it in terms of set-theoretic "correspondences" between elements. Hereafter talking about isomorphisms I shall understand this notion in the sense of the above definition.

A similar point can be made about arithmetical models of plane Euclidean and other geometries used by Hilbert in *Grundlagen*. Perhaps one can indeed imagine geometrical points as usual dots, "systems of loves" or beer mugs indiscriminately. But representation of points by pairs of real numbers (or pairs of elements of another appropriate algebraic field) is a different matter. Unlike dots and beer mugs numbers are mathematical objects on their own rights belonging to a different mathematical theory, namely arithmetic. "Translations" of geometrical theories into arithmetic used by Hilbert are obviously non-reversible: they allow for translation of geometrical theories into arithmetic terms but not the other way round. Hilbert certainly saw this. He didn't mean to say that geometry and arithmetic seen from a higher viewpoint turn to be the same theory; actually he considered a possibility of reduction of geometry to arithmetic. Nevertheless he thought about this translation as an isomorphism, namely an isomorphism between basic geometrical objects and relations, one the one hand, and specially prepared arithmetical constructions, on the other hand. But such constructions obviously cannot be made outside an appropriate arithmetical theory. As far as this "target" arithmetical theory is wholly taken into consideration (as it must be!) the translation in question doesn't look like like an isomorphism any longer. We see that the example of projective duality mentioned by Hilbert (to leave alone "systems of love, of law or of chimney sweeps") is special and cannot be used as a model for treating the notion of interpretation in mathematics in the general case. One may perhaps remark that the point I am making is obvious but trivial or at best merely technical: the categorical notion of morphism allows indeed for a better treatment of interpretations between theories but doesn't essentially change anything. Let me now show that in fact it does. What might look like a minor technical amendment suggests a revision of the whole idea of "formal" mathematics and formal theory-building in general. (Endnote 1)

## 5. Forms, Categories and Structuralism

Many traditional mathematical concepts have the following property: all items falling under a given concept are isomorphic or, in other words, are defined "up to isomorphism". I shall call such concepts form-concepts or simply forms. The colloquial expression "up to isomorphism" apparently involves a systematic ambiguity between identity and isomorphism but it doesn't matter here. Think about any Euclidean geometrical form (shape) like that of circle. A circle allows for (Euclidean) motions and scalings. I leave aside the tricky question of whether these transformations preserve identities of circles or rather associate with given circles some other circles (see Rodin, forthcoming). I only point here to the fact that all (Euclidean) circles are isomorphic in the sense that for any given pair of circles there always exist a reversible transformation (motion, scaling or their composition) transforming one circle into the other. Through an appropriate modification of the class of admissible transformations one may modify a given form-concept. Thinking about a circle up to reversible continuous transformation one gets a more general topological notion of circle. Klein (1872) first put forward the idea of description of geometrical spaces through group-theoretic properties of transformations available in these spaces. This approach requires all transformations in question to be reversible (to be isomorphisms) for otherwise they don't form a group by composition.

For a further example consider the traditional notion of (natural) number conceived as a "shared form" of isomorphic finite sets. (Number $n$ is a form shared by all sets of exactly $n$ elements; any pair of such sets allows for a one-to-one correspondence between their elements.) The colloquial notion of "algebraic form" provides yet another class of elementary examples of form-concepts. Take expression like $(x+y)\exp n$. All expressions obtained from this one through any admissible substitution of symbols $x, y, n$ by some other variables, constants or algebraic expressions are told to share the "same algebraic form". Isomorphisms associated with this algebraic form are such substitutions (notice the reversibility of substitutions). For more modern examples of form-concepts think about any concept build up as "structured set". Think, for example, about "the" group of plane Euclidean motions already mentioned. Obviously one can take as many isomorphic "copies" of this group as one likes. The mathematical notion of form (and the notion of group of symmetry associated with it) plays a major role in physics and other sciences (for a popular account see Weyl 1952).

The above examples are so various that one might think that in fact all mathematical concepts are form-concepts. This view has been held by Plato who, however, made a distinction between mathematical forms and forms strictu sensu. In todays philosophy of mathematics this view is known under the name of mathematical structuralism. Mathematical structuralism exists in a number of different versions but that structures are things determined up to isomorphism seems to be a common assumption. Understandably Hilbert is often referred to as one of founding farthers of mathematical structuralism (Hellman, forthcoming).

But the view according to which all mathematical concepts are form-concepts is obviously wrong. Think about the general concept of group - not any particular group like the group of Euclidean motions but the concept of group as such. There are certainly isomorphic groups but not all groups are isomorphic. So the general concept of group is not determined up to isomorphism. Hence it is not a form-concept. Similarly general concepts of set, natural number, etc. are not form-concepts.

This sounds like a trivial point. However one may still argue that all mathematical objects are instances of form-concepts. However important the general concept of group might be one arguably doesn't need in mathematics anything like "general group" over and above all particular groups. Since these particular groups are all form-concepts the structuralist view remains plausible. But the situation can be seen differently. Notice that like in the case of any "particular" (that is, specified up to isomorphism) group there is a notion of transformation associated with the general concept of group. I mean the notion of group homomorphism (I shall explain it shortly). Homomorphisms are non-reversible (except when they are isomorphisms!). One might say that the general concept of group is determined "up to homomorphism" like the group of plane Euclidean transformations (or any other particular group) is determined up to isomorphism. But a better way to put is this: groups with group homomorphisms make a category. The category of groups is evidently something "over and above" all particular groups that it comprises; it has properties, which cannot be detected through studying particular groups. Importantly the concept of category of groups is not a form-concept as we shall see. The same holds for the category of sets, etc. By analogy with form-concepts I shall speak about category-concepts. General concepts of sets, group, manifold, topological space or the general concept of category itself are category-concepts. This means that items falling under such concepts and their associated transformations make categories (in particular all possible categories make the category of categories, see (Lawvere 1966)). The notion of category-concept is more general than that of form-concept: any form-concept is a special case of category-concept where all morphisms are isomorphisms. (Awodey 1996) and some other people argue that the notion of category is a typical example of structure , which appears to be a different name for (or at least a special case of) what I call here form. Since the general concept of category is instantiated by sets, groups, etc., etc., an abstract category can be colloquially called a "common conceptual form" of all these things . However this is misleading and, as far as one assumes the proposed understanding of the notion of form, simply wrong. For the categories in question are not isomorph. They can be mapped to each other by suitable functors but these functors are never reversible. So it is misleading to interpret the fact that sets, groups, etc. all make categories in the sense that all these thing share a "form of category". The erroneous thinking about categories as forms is, in my view, responsible for the infamous description of category theory as "abstract nonsense".

6. Transformations Instead of Relations?

The notion of transformation can be formulated in mathematics most easily in the case when transformed objects are construed *a la* Bourbaki as "structured sets". For "bare" ("non-structured") sets we have the notion of set isomorphism as one-to-one correspondence between their elements and a more general notion of morphism between (or "transformation of") sets, which are called functions: each element of the domain is sent to a certain (unique) element of co-domain. Since two different elements of the domain can be sent to the same element of the codomain morphisms (transformations) of sets are, generally speaking, non-reversible. When a given set is equipped with a certain structure making it into a group or something else the permissible transformations are those transformations of underlying sets (that is, functions), which are said to "preserve" the corresponding structure. For example the group structure on a given set is a binary operation defined for elements of this sets and verifying certain axioms. Given two groups *G* and *G'* construed as just explained with corresponding group operations + and * think about function *f*: *G*-->*G'*. It is said to preserve the group structure (and hence to be a group homomorphism) iff

(GS):   $f(x+y)=f(x)*f(y)$.

When *f* is reversible groups *G* and *G'* are called isomorphic or "identical up to isomorphism". Such groups cannot be distinguished by their group-theoretic properties. However (GS) may hold in the non-reversible case as well. Then the colloquial talk of "preservation of structure" becomes rather misleading. Consider the case when *G'* consists of only one element 1 (such that 1+1=1). In this case homomorphism *f* "destroys" the structure of *G* rather than preserves it! A better way to put (GS) into words is to say that in the general case homomorphisms *respect* the group structure (even if they destroy it).

Things work similarly with all "structured sets" (with different conditions amounting to "preservation" or "respect" of the corresponding structure). This points to an important link between category theory and the notion of mathematical structure (and hence with Mathematical Structuralism). However I see this link as historical rather than theoretic. For the notion of transformation in fact doesn't depend of the set-theoretic background involved in definitions of morphisms of groups, topological spaces and other concepts construed as structured sets. This is true in the case when transformations are conceived "naively" as well as in the case when the notion of transformation acquires a technical definition other than set-theoretic. The notion of geometrical transformation came about long before set theory; in particular the notion of geometrical motion is implicit already in Euclid's *Elements* (think about Euclid's notion of congruence). Klein in his *Erlangen Program* made geometrical transformations into foundations and treated them algebraically using the notion of (algebraic) group. In this latter context the notion of transformation could be hardly called naive any longer. Finally the notion of transformation or morphism became a basic notion of category theory. Although in the early days of category theory all important examples of morphisms were morphisms of structured sets MacLane and Eilenberg realised it from the very beginning that certain morphisms are not

"structural" in this sense. Two basic examples of "non-structural" categories which immediately suggest themselves are the following: (1) group construed as a category with just one object and such that all its morphisms are isomorphisms; (2) preoder construed as a category having no more than one morphism between any (ordered) pair of its objects. The categories just mentioned shouldn't be confused with categories of "all" groups and "all" preorders correspondingly since (1) is an "individual" group and (2) is an "individual" preorders. But category of "all" groups and category of "all" preorders could be construed by taking, correspondingly, (1) and (2) as their objects. Groups and preorders, of course, "are" structured sets in the sense that they can be construed as such things. However as we have just seen they can be also construed "purely categorically", so their morphisms are not supposed to "preserve" or "respect" any structure. However trivial may be the translation of notions of group and preorder from the set-theoretic into the categorical language it clearly shows that the notion of morphism is more general than that of structure-preserving morphism.

Let's now go back to formal axiomatic method. Its basic idea can be expressed by this slogan: describe mathematical objects in terms of formal properties of their relations. The slogan of categorical method reads slightly differently: describe mathematical objects in terms of categorical properties of their transformations. The principle assumption behind the categorical approach is that transformations of mathematical objects (including both isomorphisms and non-reversible transformations) indeed essentially characterise these objects.

A comparison between the two methods suggests an analogy between relations and morphisms. One can say indeed that morphisms "relate" objects to each other in a way. As far as this claim is taken in the general philosophical sense it sounds reasonable. But it is in odds with the standard technical notion of relation as function sending tuples of individuals (relata) to truth values. For morphism $f:A \rightarrow B$ sending object $A$ to object $B$ is just another mathematical object, which prima facie has nothing to do with truth values. Further, in a given binary relation $R(x,y)$ arguments (relata) $x,y$ can be, generally speaking, replaced by some other arguments $x',y'$ so that $R(x',y')$ again makes sense (i.e. is true or false). This allows for thinking about $R$ as "structural relationships" (between its relata) appearing in Hellmann's official definition of structuralism as a view about the subject matter of mathematics according to which what matters are structural relationships in abstraction from the intrinsic nature of the related objects

But for morphism *f:A-->B* nothing similar is possible. The fact that object *A* is domain of *f* and object *B* is codomain of *f* characterises *f* essentially: there is no sense in which *f* may survive replacement of *A*, *B* by some other objects *A'*, *B'*.
These observations point to a gap between the general philosophical intuition behind the concept of relation and the standard logical notion of relation just mentioned. The intuition says that morphisms are relation-like while within the standard logical formalism they should be treated as objects. This can be seen in Lawvere's early papers of 1964 and 1966 already mentioned where categories are treated as classes of morphisms with relations "domain of", "codomain of" and "composition of". Let me now provide a further argument showing that doing category theory formally (i.e. using formal axiomatic method) is not a good idea.

7. Categories and Categoricity.
As the above quote from Hilbert clearly shows he thinks of a formal theory as construed up to isomorphism of its interpretations. But what guaranties that a given theory like that of *Grundlagen* is indeed formal in this sense, i.e. that all its models are in fact isomorphic? The desired property has been called by Veblen (1904) categoricity. This term has nothing to do with category theory.
When Hilbert was preparing his *Grundlagen* for publication he apparently didn't yet see the problem. He discovered it about the time of the first publication of *Grundlagen*. In his lecture *Ueber den Zahlbegriff* delivered in 1899 and published in 1900 Hilbert first introduced an "axiom of completeness" (Vollstandigkeitsaxiom) requiring from any model of a given theory (this time it was arithmetic) this maximal property: given model *M* satisfying the rest of the axioms one cannot obtain another model satisfying the same axioms by extending *M* with new elements. In the second and following editions of *Grundlagen* Hilbert used a similar axiom under the same name. Hilbert's Vollstandigkeitsaxiom implies categoricity (although the converse is obviously not the case). From the point of view of today's model theory this axiom looks very dubios if not plainly "wrong": Hilbert's account doesn't provide any reason why a model with the desired property should exist (in any appropriate sense of "exist") but apparently relies onto the intuition which suggests that the intended model (i.e. "usual" geometrical space) has this property. This is a very shaky ground indeed. The standard Tarskian model theory doesn't allow for a model with the required maximal property (because of "upward" Skolem-Lowenheim's theorem: given that a theory has an infinite model it has other models in higher cardinalities). So the Vollstandigkeitsaxiom turned to be incompatible with later model theory. Since the notion of categoricity has been formulated by Veblen the idea that this property can be stipulated by a fiat has been largely abandoned.

For the reason just explained categoricity is commonly viewed as a desired property of formal theories. However in 20-th century people have learnt to be tolerant to the lack of categoricity. For Zermelo-Frenskel axiomatic set theory (ZF), Peano Arithmetic (PA) and some other theories commonly viewed as important turned to be non-categorical. To preclude the right of these theories to be qualified as formal on this ground would apparently mean to go too far. To save the situation philosophers invented the notion of "intended model", that is of model chosen among others on an intuitive basis. Isn't this ironic that such a blunt appeal to intuition is made in the core of formal axiomatic method? I'm agree with F. Davey who recently argued that "no-one has ever been able to explain exactly what they mean by intended model". (FOM, 13 Oct 2006). Other people question the categoricity requirement. Asks R. Lindauer: "Why rule out non-standard models of 1st-Order PA? What's wrong with having other models? Why should we be making our model-world smaller and not larger?" (FOM, 27 Oct 2006).

I believe that the lack of categoricity of theories like ZF and PA is indeed a serious flaw because the lack of categoricity undermines the very idea of formal theory. At the same time I agree with Lindauer and other people who think that the pursuit of categoricity is misleading. These two claims might seem to contradict each other but they don't: instead of forcing categoricity or looking for a philosophical excuse of the lack of categoricity of formal theories I suggest to change the method of theory-building and the corresponding notion of theory. As far as non-reversible morphisms are treated on the same footing with isomorphisms the pursuit of categoricity has no sense any longer. Trying to describe a model of a given theory "up to arbitrary morphism" rather than up to isomorphism one may get a category of models which has

"good" categorical properties making it "well-manageable". In the following paragraph I provide some details of how this can be achieved. As we shall see the categorical approach undermines the usual distinction between a (formal) theory and its models: in the new context a theory can be naturally seen as one of the models having this specific property that it "generates" all the others. This is hardly surprising given that the notion of formal theory (as distinguished from its models) requires categoricity in Veblen's sense.

A great advantage of formal axiomatic method is that it provides a clear idea about the role and the place of logic of in theories built by this method (albeit details can always be a subject of philosophical discussion). Talking about categorical method as an alternative to formal method we cannot avoid this important issue either. In the next paragraph I shall briefly review classical and formal notions of logic, and then develop a notion of categorical logic and show its role in categorical method of theory-building.

## 8. Formalising Logic

Traditionally logic is conceived as a general theory of reasoning independent of any particular subject matter one may reason about. On Aristotle's account logic is closely related to ontology. In particular, Aristotle treats the logical law of non-contradiction as a fundamental ontological principle, and his so-called perfect syllogism reflects the structure of entity (as Aristotle understands it). Logical truths are grounded upon ontological truths even if the former do not coincide with the latter. There is, of course, a sense in which Aristotle's logic could be called formal. For it captures and studies common forms shared by various reasonings about different matters. These forms are called logical forms; perfect syllogism is a typical example. However other sciences like biology proceed similarly: biology captures forms shared by different organisms and so brings about the notion of biological form (differently called "form of life" or "living form"). However biology can be hardly called a formal science on this ground. This shows that the notion of being formal relevant to modern logic is different.

The Hilbertian notion of formal theory is that of "framework or schema of concepts together with their necessary relations to one another" taken in abstraction from its possible "basic elements", that is, technically speaking, from its possible interpretations. Formal logic in the usual today's understanding of this expression is formal in the same sense. This gives the idea of distinction between logical syntax and logical semantics which is not found in the traditional logic. "Formal" means here "syntactic". As Carnap puts this

The task of formalisation of any theory ... belongs to syntax, not to semantics. (1947, Preface)

However there is a problem here, which make it difficult to apply the notion of formal theory to logic. Remind that not all terms used in axioms of Hibert's *Grundlagen* have variable meanings. Meanings of terms "and", "or", "exist" and of some others are fixed; such terms form the logical vocabulary of the given theory, and the (maximal) fragment of the theory which involves no other terms but logical can be identified with logic. So unlike geometrical theories themselves their underlying logic is fixed and doesn't allow for different interpretations. But this means that one cannot distinguish here between formal and interpreted logic or between logical syntax and logical semantics. However with this latter distinction we loose the modern notion of formal logic.

Actually this notion is apparently absent in Hilbert. The idea behind his *Grundlagen* is to base mathematics in general and geometry in particular on logical rather than intuitive or other grounds. Discovery of non-Euclidean geometries and some other developments in mathematics led many people to the belief that since the traditional mathematical intuition proved unreliable logic remained the only firm foundations available for mathematics. (Nagel 1939) provides a thorough historical analysis of mathematics of 19th century showing where this logicist view on mathematics stems from (Endnote 2). Hilbert's *Grundlagen* showed how this general approach can be realised in practice. A formal theory in Hilbert's sense is a "logical skeleton" (or "logical form" liberally understood) shared by a class of traditional so-called "naive" mathematical theories. But logic itself on this account is not formal in anything like the same sense. Applying the notion of formal theory just given to logic one would need to speak of "logical skeleton of logic" which is at least unclear and at most senseless.

Carnap and other promoters of the idea of formal (or "formalised" logic) largely disregarded this philosophical difficulty and applied formal or "syntactic" method to logic itself, introducing the nowadays standard distinction between logical syntax and logical semantics. Instead of saying that in formal theories of the type of *Grundlagen* logic is not a subject of interpretation these people (as well as many today's professors of logic after them) would say that logical terms get interpreted together with non-logical terms but unlike the latter they are always interpreted in the same way (i.e. logical terms are invariant under all possible interpretations). Unless the class of "all interpretations" is precisely determined this move seems me purely rhetorical, and if such a class is determined then the notion of logicicity becomes relational (dependent of the given class of interpretations). Anyway this doesn't solve the problem, which is the following.

The initial hope that logic unlike geometry will always be rigidly fixed on the pain of absurdity turned to be futile and logic ramified into multiple systems just like earlier did geometry. Formalisation of logic played an important role in this development because it allowed to treat systems of logic on equal footing with systems of geometry or algebraic systems. However since the assumption about rigidity of logic is given up the whole idea of formal approach (at least in Hilbert's sense) is shaken, so it becomes rather unclear what is meant by "formal" logic except that this kind of logic is symbolic and mathematical. What kind of new philosophy of logic is needed to replace Hilbert's (or Frege's) logicism in order to cope with these developments remains an open question. Most philosophers working today in logic share Hilbert's weak logicism according to which logic has to do with foundations of mathematics and of other sciences. At the same time only few of them if any hold the old-fashioned view according to which there is only one "true" system of logic. As the above analysis shows the two parts of the popular view are hardly compatible with each other. Categorical logic suggests a solution of this problem through a revision of formal axiomatic method and more broadly - of formal approach. In order to show this solution I shall first develop a speculative notion of categorical logic as a generalisation of formal logic, and then point to some technical developments supporting this speculative notion.

## 9. Categorical Logic

Both traditional Aristotelian logic and modern formal logic hinge on the notion of logical form. What kind of forms are logical forms is a difficult question which I shall not now try to answer here. Let's see instead what happens to logic when the notion of form is upgraded to that of category. Remind that categories unlike forms, generally speaking, don't allow for a straightforward abstraction: given a class of balls one may think about them "up to isomorphism" and stipulate The Ball as their shared abstract form but nothing similar works when objects of a given class make a category. So a categorical system of logic unlike formal logic cannot be anything like a self-standing structure occasionally applied in this or that particular context. Instead it must be "internal" or "intrinsic" with respect to a given category playing the role of such context. This rises the question of universality of categorical logic: Is this indeed appropriate to give the title of logic to something, which applies to a particular category rather than to everything? Let me make three remarks concerning this question. First, nothing prevents one to conceive of "everything" as a category (rather than as a class). This idea is behind (Lawvere 1966) "The Category of Categories" approach. Personally I'm not sympathetic with this idea. Actually I consider the "local" character of categorical logic as its advantage rather than otherwise. The second remark is that the idea of "local" or "regional" logic as a notion of a system of logic designed for some specific purposes has been around already during few decades, and it better fits today's technical developments in logic than the traditional idea of the universal logic. The third remark is that in the categorical setting the notion of "regional" or "local" logic can better cope with the following important objection: Ramification of logic brakes the rational thought into a number of incompatible domains and this contradicts the whole idea of rationality. The usual response to this problem is the suggestion to find a weak system $U$ of universal logic such that regional logics could be seen as specifications of $U$ in corresponding local contexts. Alternatively one may challenge the assumption about the unity of rationality on philosophical grounds. Categorical logic allows for a different solution, namely it provides means for translation between different local logics. Such translation doesn't intervene here as a new external principle since what I call here local or regional logic is construed in the categorical setting in terms of morphisms which can be naturally viewed as translations. One would still need, of course, some universal principles, namely the general principles of category theory. What makes the difference is this: in the categorical setting universal principles of rationality are principles of translation rather than formal principles imposing universal forms of reasoning indifferent to its content. One may argue that general principles of translation I'm talking about are themselves formal but this is, in my view, an abuse of the language. As far as one tries to be precise about the meaning of "formal" it becomes clear that the argument is wrong. Importantly categorical logic assumes the possibility of multiple local logics to begin with, so that no counterpart of the aforementioned problem about formalisation of logic arises in the categorical context.

Let's now see how this speculative notion of categorical logic can be realised technically. There are several different ways to "do logic" with categories but the most relevant in the present context is apparently the so-called topos logic. The notion of topos is of a geometrical origin; it was a discovery of Lawvere that this notion can be introduced axiomatically through an appropriate specification of the abstract notion of category. It turns out that given a topos one may associate with it a logical calculus called "internal language". Then the given topos can be viewed as a geometric model of this calculus. However a more suggestive view on this situation is the converse one: the given topos has certain specific "logical properties" which determine its "internal logic".  This gives the idea of "reasoning in a topos"; reasonings in different toposes can be always compared through morphisms (functors) between these toposes (which can be of different sorts). Remark that the view just mentioned (colloquially known as "toposophy") is in odds with the usual (weak) logicism which requires to "fix logic first". For according to toposophy logic is an element (or perhaps an aspect) of the overall construction of topos, which doesn't have an epistemic priority with respect to the rest of this construction.  This rises anew the traditional philosophical issue about first principles. But let me now turn this discussion in a more technical mode.

Technically speaking the problem is that the general idea of categorical logic doesn't provide by itself any clue of  how it could be used for theory-building. Even if one refuses the idea of logical foundations of theories and tries to recover logic afterwards the problem of foundations still persists at least as a pedagogical problem. For no theory can be grasped at once but needs certain guiding mechanisms allowing to explore it piece by piece. In addition any theory needs an entry (or multiple entries). Let me briefly describe the idea of functorial semantics put forward in (Lawvere 1963), which serves these purposes. It can be viewed as a simulation of formal axiomatic method by categorical means. This view on functorial semantics is helpful for comparing formal and categorical methods.

Instead of writing axioms with usual strings of formulas one encodes axioms into a special "syntactic" category **T** which plays the role of "formal theory". Then like in the case of standard semantics one takes for granted a "background category" **B**, which is usually taken to be the category of sets but can be chosen differently. Models of **T** are functors of the form **F**: **T**-->**B**. This construction is that it allows for different notions of model dependent one specific properties of functors of form **F**. Another remarkable fact is that under rather general conditions **T** can be embedded into a category **M**(**T**,**B**) of its functorial models. This definitely changes the whole idea of theory as a structure over and above all its possible models and suggests the view on a theory as "generic model" (Lawvere 1963-2004, p.19) which generates other models like circles and straight lines generate further constructions in Euclid's *Elements*. The functorial semantics makes it clear that the requirement of categoricity (in Veblen's sense) is as much unrealistic as unnecessary: although "good" categorical properties of **M**(**T**,**B**) are much desirable there is no good reason to require that this category reduces to a single object.

10. Conclusion: Formalisation versus Categorification

Lawvere's functorial semantics has been developed for a special case of algebraic theories and so it cannot be immediately used as a method of theory-building applicable in all areas of mathematics. Since then a lot of technical work been done in related fields of categorical logic and categorical model theory. For a historical introduction and further references I refer the reader to (Bell 2005). The categorical method of theory-building is a work in progress, so to the date it doesn't exist in any standard form. The purpose of this non-technical paper is to provide this method with an appropriate epistemological background, which might motivate further technical work. So let me conclude with some general epistemological remarks.

The categorical method outlined above suggests an epistemic strategy, which differs from that suggested by formal method. In the latter case the general epistemic strategy is to subsume different objects under a common form and then stipulate this form as a self-standing abstract object. This amounts to identifying of given objects "up to isomorphism". Here "objects" may stand for various mathematical constructions including whole mathematical theories, so the sense in which the obtained formal object is "abstract" is relational. This leads to a traditional hierarchical structure of organisation of mathematical knowledge, where theories and concepts are subsumed under other theories and concepts, which are more abstract, more general and more formal. Apparently the same pattern applies at the interface between mathematics and the material world and is responsible for the usual qualification of mathematics as a formal science. This epistemic strategy I call formalisation.

The alternative strategy of categorification is different. It is more general and in a sense more straightforward. Categorification amounts to taking into consideration all transformations, which can be described as categorical morphisms but not only reversible ones, i.e. not only isomorphisms. The requirement that objects and transformations in question make a category is much weaker than the requirement that these transformations are reversible (and so make a group). So formalisation is a special case of categorification. Unlike formalisation categorification in the general case doesn't subsume objects in question under a conceptual umbrella but simply links them (by mophisms) into into a whole, namely into a category. When objects are theories and morphisms are mutual interpretations of theories one gets a network of theories. Although this network might have no single "centre" it may be still coherent and well-manageable if it has "good" categorical properties.

We see that categorification like formalisation serves for integration of mathematical knowledge. But categorification unlike formalisation doesn't bring about a hierarchical structure. Given that the very idea of foundations seems to imply a hierarchical structure of knowldge (which "starts with" foundations and then branches into various specific sub-domains) one may wonder if categorification is compatible with it. I think that at least one version of the notion of foundations remains viable in a categorical context, namely the pedagogical one stressed by Lawvere in his (2003). I mean the notion of foundations as an "entry" into a theoretical network.  Such entry should exist for any theoretical network but it obviously needs not to be unique.

I suggest that hierarchical structures cannot any longer serve as universal models of organisation of knowledge just like they cannot any longer serve as models of organisation of our societies. But the task of integration of knowledge into a manageable whole remains pertinent as ever. I believe that categorical method of theory building outlined in this paper can be helpful for this task.

Endnotes:

(1) One might argue that I systematically confuse here two different notions of interpretation: (a) assignment of referents to primitive terms of a formal theory like "point" and "straight line" (evaluation of logical variables) and (b) interpretation of one non-formal theory in terms of another non-formal theory like in the case when traditional geometrical points are represented by pairs of numbers. In fact in preceding paragraphs I'm talking only about (b) as also does Hilbert in the given quote. The argument which I develop below in the main text is this: the idea of formal theory and of interpretation in the sense (a) assumes that all models of a given formal theory are mutually reversibly interpretable in the sense (b). But generally they are not.

(2) This logicism about mathematics should be distinguished from stronger form of logicism aiming at reduction of mathematics to logic.

Literuature: